\documentclass{amsart}

\usepackage{amssymb}
\usepackage{latexsym}
\usepackage{xypic}

\DeclareMathAlphabet{\mathrc}{U}{eur}{m}{n}

\newcommand{\GL}{\operatorname{GL}}

\newcommand{\Hom}{\operatorname{Hom}}
\newcommand{\End}{\operatorname{End}}
\newcommand{\tensor}{\otimes}
\newcommand{\soc}{\operatorname{soc}}

\newcommand{\res}{\operatorname{res}}

\newcommand{\oo}{\mathfrak{o}}  
\newcommand\R{\mathbb{R}}    
\newcommand\Z{\mathbb{Z}}    
\newcommand\F{\mathbb{F}}    
\newcommand\Q{\mathbb{Q}}    

\newcommand\Lie{\operatorname{Lie}}
\newcommand\Ad{\operatorname{Ad}}

\newcommand\lie[1]{\mathfrak{#1}}
\newcommand\glie{\lie{g}}
\newcommand\ulie{\lie{u}}

\newcommand{\NN}{\mathcal{N}}
\newcommand{\normal}{\lhd}
 
\newcommand{\iso}{\simeq}
\newcommand{\GG}{\mathbf{G}}
\newcommand{\da}{\bullet}

\swapnumbers

\theoremstyle{plain}

\newtheorem*{theorem}{Theorem}

\swapnumbers
\newtheorem{theoreml}{Theorem}

\swapnumbers

\newtheorem*{prop}{Proposition}

\newtheorem*{lem}{Lemma}

\swapnumbers
\newtheorem{leml}{Lemma}

\swapnumbers

\swapnumbers

\swapnumbers
\theoremstyle{remark}
\newtheorem*{rem}{Remark}
\newtheorem*{rems}{Remarks}

\numberwithin{equation}{subsection}

\begin{document}
\bibliographystyle{amsalpha} 
\author{George J. McNinch}
\thanks{This work was supported by a grant from the 
        National Science Foundation.}
\email{McNinch.1@nd.edu} 
\date{\today}

\begin{abstract}
  Let $G$ be a simple, simply connected and connected algebraic group
  over an algebraically closed field of characteristic $p>0$, and let
  $V$ be a rational $G$-module such that $\dim V \le p$.  According to
  a result of Jantzen, $V$ is completely reducible, and $H^1(G,V)=0$.
  In this paper we show that $H^2(G,V) = 0$ unless some composition
  factor of $V$ is a non-trivial Frobenius twist of the adjoint
  representation of $G$.
\end{abstract}

\renewcommand{\subjclassname}{%
    \textup{2000} Mathematics Subject Classification}

\address{Department of Mathematics \\ 
  University of Notre Dame \\ Notre Dame, IN~ 46556 ~USA}

\title[The second cohomology of small modules]
      {The second cohomology of small irreducible modules for
        simple algebraic groups}

\maketitle

\section{Introduction}

Let $G$ be a quasisimple, connected, and simply connected algebraic
group over the algebraically closed field $k$ of characteristic $p>0$.
By a $G$-module $V$, we always understand a rational $G$-module (one
given by a morphism of algebraic groups $G \to \GL(V)$).  In this
paper, we study the cohomology of a $G$-module $V$ such that $\dim V
\le p$. By results of Jantzen \cite{J3a} one knows that $V$ is
semisimple and that $H^1(G,V) = 0$.

Recall that the Lie algebra $\glie$ of $G$ is a $G$-module
via the adjoint action.  Our main result is:
\begin{theoreml} 
  \label{intro:main-theorem}
  Let $V$ be a $G$-module with $\dim V \le p.$ Then $H^2(G,V) \not =
  0$ if and only if $V$ has a composition factor isomorphic with a
  \emph{Frobenius twist} $\glie^{[d]}$ of $\glie$ for some $d \ge 1$.
\end{theoreml}
  
Differentiating the representation of $G$ on $V$ gives a
representation for the Lie algebra $\glie$ on $V$. Assume that
$V^\glie = 0$.  Then the theorem says that $H^2(G,V) = 0$. For $V$ of
this sort, the vanishing of $H^2$ is a consequence of the linkage
principle for $G$ together with results in section \ref{sec:root-sys}
which give estimates for the dimensions of Weyl modules whose high
weights are simultaneously in the low alcove and in the orbit $W_p \da
0$. In fact, the same argument shows that $H^i(G,V)$ is 0 for all $i
\ge 1$; see Proposition \ref{sub:G-vanish-for-G1-non-triv}. It was
pointed out to me that an earlier version of this manuscript contained
an overly complicated proof of this observation.

The crucial case for Theorem \ref{intro:main-theorem} is when $V$ is
simple, non-trivial and $V^\glie = V$.  There is a unique $d \ge 1$
such that the ``Frobenius untwist'' $V^{[-d]}$ is a $G$-module on
which $\glie$ acts non-trivially.  We have already seen that
$H^i(G,V^{[-d]}) = 0$ for $i=1,2$, so Theorem A follows from the
following two results (see \ref{sub:second-vanishing}). [We denote by
$h$ the Coxeter number of the group $G$.]

\begin{theoreml}
  \label{intro:tool-theorem}
  Suppose that $p \ge h$ and that $W$ is a $G$-module for which
  $H^i(G,W) = 0$ for $i=1,2$. Then $H^2(G,W^{[d]}) \iso
  \Hom_G(\glie,W)$ for all $d \ge 1$.
\end{theoreml}
\begin{theoreml}
  \label{intro:lie-alg}
  If $p>h$, $\dim H^2(G,\glie^{[d]}) = 1$ for all $d \ge 1$.
  For any $p$, there is a $d_0 \ge 1$ so that $H^2(G,\glie^{[d]}) \not = 0$
  for all $d \ge d_0$.
\end{theoreml}
Theorem \ref{intro:tool-theorem} is proved in \ref{sub:H2-tool}; it
immediately implies the first assertion of Theorem \ref{intro:lie-alg}
(see \ref{sub:adjoint-second}). We give a proof the second assertion
of Theorem \ref{intro:lie-alg} in section
\ref{sub:finite-obstruction}.

We end the paper by applying the results of section \ref{sec:root-sys}
to calculations of cohomology groups $H^i(G_1,L)$, where $G_1$ is the
Frobenius kernel, and $L$ is a simple $G_1$ module with $\dim L \le
p$; see Proposition \ref{sec:G1}.

We conclude this introduction by remarking that the result of Jantzen
\cite{J3a} cited above is one of several recent results studying the
semisimplicity of low dimensional representations of groups in
characteristic $p$.  See \cite{serre}, \cite{M2},
\cite{M3},\cite{Gural-ss}, and \cite{Mc:SemiExterior} for related
work.

The author would like to acknowledge the hospitality of Bob Guralnick
and the University of Southern California during a visit in June 1999;
in particular, questions of Guralnick encouraged the author to
consider the problems addressed in this paper, and several
conversations inspired some useful ideas. 

\section{Root systems}
\label{sec:root-sys}

\subsection{}
We denote by $R$ an indecomposable root system in its weight lattice
$X$ with simple roots $S \subset R^+$. For each $\alpha \in S$,
there is a fundamental dominant weight $\varpi_\alpha \in X$; the
$\varpi_\alpha$ form a $\Z$ basis of $X$.  

We write $\alpha_0$ for the dominant short root, and $\tilde \alpha$ for
the dominant long root in $R$ (these coincide in case there is only
one root length).

The Coxeter number of $R$ is given by
\begin{equation*}
  h - 1 = \sup_{\alpha \in R^+}\{\langle \rho,\alpha^\vee\rangle \}
  = \langle \rho,\alpha_0^\vee \rangle.
\end{equation*}

For $m \in \Z $ and $\alpha \in R$, let $s_{\alpha,m}$ denote the
affine reflection of $X_\R = X \tensor_\Z \R$ in the hyperplane
$H_{\alpha,m} = \{ x \in X_\R: \langle x,\alpha^\vee \rangle = m\}$.

Let $l > h$ be an integer.  The affine Weyl group $W_l$ is the group
of affine transformations of $X_\R$ generated by all $s_{\alpha,ln}$
for $n \in \Z$.  According to \cite[ch. VI, \S 2.1, Prop.
1]{BouLie456} $W_l$ is isomorphic to the semidirect product of $W$
(the finite Weyl group) with $l\Z R$.  The normalizer of $W_l$ in the
full affine transformation group of $X_\R$ contains all translations
by $lX$, hence $W_l$ is a normal subgroup of $\widehat W_l$, the
semidirect product of $W$ and $lX$. Moreover, $\widehat W_l / W_l \iso
lX/l\Z R \iso X/\Z R$ is the fundamental group of $R$, which we will
denote by $\pi$.

Let $\rho = \dfrac{1}{2}\sum_{\alpha \in S} \alpha$.  We always
consider the dot action of $\widehat W_l$ (also of $W$ and $W_l$) on
$X$: for $w \in \widehat W_l$ and $\lambda \in X$, this is given by $w
\da \lambda = w(\lambda + \rho) - \rho$.

The subset $C_l$ of $X_\R$ given by
\begin{equation*}
  C_l = \{ \lambda \in X_\R \mid 
  0 < \langle \lambda + \rho, \alpha^\vee \rangle < l 
  \quad \text{for each $\alpha \in R^+$}\}.     
\end{equation*}
is a fundamental domain for the dot action of $W_l$ on $X$; its
conjugates under $W_l$ are known as alcoves, and $C_l$ is the lowest
alcove. Since $\widehat W_l$ normalizes $W_l$, \cite[ch. VI, \S
2.1]{BouLie456} shows that $\widehat W_l$ permutes the alcoves.

Let $\Omega$ be the stabilizer in $\widehat W_l$ of $C$. Since $W_l$
permutes the alcoves simply transitively, one deduces that $\widehat
W_l$ is the semidirect product of $\Omega$ and $W_l$.  Thus $\Omega \iso
\widehat W_l/W_l \iso \pi$.

Since $l > h$, the intersection $C_l \cap X^+$ is non-empty.  [Note
that if $l \le h$ had been allowed, we would have $C_l \cap X^+ =
\{0\}$ in case $l=h$, and $C_l \cap X^+ = \emptyset$ if $l < h$.]  It
is then clear that $\widehat W_l \da 0 \cap C_l = \{\omega \da 0 \mid
\omega \in \Omega\}.$

\subsection{} 
\label{sub:fund-group}

Let $I$ index the simple roots $S = \{\alpha_i\}$, write
$\alpha_0^\vee = \sum_{i \in I} n_i \alpha_i^\vee$, and put $J = \{ i
\in I \mid n_i = 1\}$.  A dominant weight $0 \not = \varpi \in X$ is
\emph{minuscule} if whenever $\lambda \le \varpi$ and $\lambda$ is a
dominant weight, then $\varpi = \lambda$. According to \cite[Ch. VI,
exerc.  23,24]{BouLie456}, $\varpi$ is minuscule just in case $\varpi
= \varpi_i$ for some $i \in J$.

For $i \in I \cup \{0\}$, let $S_i = S \setminus \{\alpha_i\}$ (so
$S_0 = S$).  Write $R_i \subset R$ for the root subsystem determined
by $S_i$, and $W_i$ for the parabolic subgroup of $W$ associated with
$R_i$.  Let $w_i \in W_i$ be the unique element which makes all
positive roots in $R_i$ negative.

For $x \in X$, let $t(x)$ denote the affine translation by $x$; for $i
\in J$, let $\gamma_i = t(l\varpi_i)w_0w_i \in \widehat{W_l}$.  Note
that $\gamma_i$ represents $\varpi_i \in X/\Z R \iso lX/l\Z R \iso
\widehat{W_l}/W_l.$

Applying \cite[ch. VI, \S2.2 Prop.  6 and Cor.]{BouLie456} one obtains:
\begin{prop}
  \begin{enumerate}
  \item[(a)] Each non-0 coset of $\Z R$ in $X$ is uniquely represented
    by a minuscule weight. In particular, $|\pi| = |J|+1$.
  \item[(c)] The non-identity elements of $\Omega$ are precisely
    the $\gamma_i$ for $i \in J$. We have
    \begin{equation*}
      \widehat W_l \da 0 \cap C_l = 
      \{0\} \cup \{\gamma_i \da 0 = (l-h)\varpi_i \mid i \in J\}
    \end{equation*}
  \end{enumerate}
\end{prop}

\subsection{}
\label{sub:Weyl-degree-estimates}
For a dominant weight $\lambda$, let
\begin{equation}
  \label{eq:Weyl-degree}
  d(\lambda) = \prod_{\alpha > 0} \dfrac{\langle \lambda + \rho, \alpha^\vee
\rangle}{\langle \rho, \alpha^\vee \rangle}
\end{equation}
be the value of 
Weyl's degree formula at $\lambda$.

\begin{prop}
  \label{prop:degree-estimates}
  Let $\lambda = (l-h)\varpi_i$ for some $i \in J$.
  \begin{itemize}
  \item[(a)] $d(\lambda) \ge \binom{l-1}{l-h}$, with equality if and only if
    $h -1 = \ell(w_0w_i)$.
  \item[ (b)] If $l-h\ge 2$ and $h \ge 3$, then $d(\lambda) > l$.
  \end{itemize}
\end{prop}

\begin{proof}
  For $1 \le k \le h-1$, let $e(k)$ be the number of $\alpha \in R^+
  \setminus R_i^+$ with $\langle \rho, \alpha^\vee \rangle = k$.  The
  argument in the remark on p. 520-521 of \cite{serre} (following
  Prop. 6) shows that $e(k) \ge 1$ for each $1 \le k \le h-1$.  Thus,
  we have
  \begin{equation*}
    d(\lambda)  = \prod_{k=1}^{h-1} \left( \dfrac{l-h+k}{k} \right)^{e(k)}
    \ge \prod_{k=1}^{h-1} \dfrac{l-h+k}{k} = \dbinom{l-1}{l-h}.
  \end{equation*}
  If $\ell(w_0w_i) = |R^+| - |R_i^+| = h-1$, then $e(k)=1$ for each $1
  \le k \le h-1$ and equality holds. This proves (a).
  
  For (b), note that under the given
  hypothesis we have $l \ge 5$. Since $\binom{l-1}{l-h} \ge
  \binom{l-1}{2} > l$ for all such $l$, (b) follows immediately.
\end{proof}

\begin{rem}
  Using the table in the proof of Proposition
  \ref{sub:Weyl-degree-exceptions} below, it is straightforward to
  verify that equality holds in (a) if and only if either $R=A_r$ and
  $i \in \{1,r\}$ or $R=C_r$ and $i = 1$. (Since $B_2 = C_2$, the
  latter case includes $B_2$ and $i=2$.)
\end{rem}

\subsection{}
\label{sub:Weyl-degree-exceptions}

In the following, let me emphasize  the standing assumption $l > h$.
\begin{prop}
  If $0 \not = \lambda \in \widehat W_l \da 0 \cap C$ and $d(\lambda)
  < l$ then $d(\lambda) = \ell -1$ and $(R,\lambda)$ is listed in the
  following table. If the rank of $R$ is $\ge 2$, then $l = h+1$.
\end{prop}

\begin{equation*}
  \label{eq:small-possibilities}
  \begin{array}[t]{lll}
    R & l & \lambda \\
    \hline \hline
    A_1 & \text{any} & (l-2)\varpi_1 \\
    A_{l-2} & & \varpi_1, \varpi_{l-2} \\
    B_2 & l = 5 & \varpi_2 \\
    C_{(l-1)/2} & \text{$l$ odd}  & \varpi_1 
  \end{array}
\end{equation*}

\begin{proof}
  The rank 1 situation leads to the item listed in the table. When the
  rank is at least 2, one applies Proposition \ref{prop:degree-estimates} to
  obtain $l = h+1$, whence $\lambda = \varpi_i$ for some
  $i \in J$; i.e.  $\lambda$ is minuscule.
  
  We handle the minuscule cases by classification.  For each
  indecomposable root system $R$ for which $J \not = \emptyset$, we
  list in the following table the Coxeter number, the set $J$, and the
  value $d(\varpi_i)$ for each $i \in J$.  The simple roots are
  indexed as in the tables in \cite[Planche I-X]{BouLie456}; the data
  recorded here, with the exception of the values $d(\varpi_i)$, may
  be verified by inspecting those tables as well.  The values 
  $d(\varpi_i)$ are well known (and can anyway be computed from the
  formula, or by representation theoretic arguments).
  \begin{equation*}
    \begin{array}[t]{l|llll}
      \text{Type of $R$} & h & J & d(\varpi_i), \ i \in J \\
      \hline \hline
      A_r & r+1 & \{1,2,\dots,r\} & \binom{r+1}{i} \\
      B_r, r \ge 2 & 2r &  \{r\} & 2^r \\
      C_r, r \ge 2 & 2r  & \{1\} & 2r \\
      D_r, r \ge 4 & 2r-2 & \{1,r-1,r\} & 
                            2r,2^{r-1},2^{r-1} \ \text{respectively} \\
      E_6 & 12 & \{1,6\} & 27,27 \\
      E_7 & 18 & \{7\} & 56  
    \end{array}
  \end{equation*}
  From this table, one can list all pairs $(R,\lambda)$ for which $R$
  has Coxeter number $l-1$ and $\lambda$ is minuscule. It is a simple
  matter to see that $d(\lambda) < l$ only when $(R,\lambda)$ is as
  claimed.
\end{proof}

\section{The algebraic groups}

\subsection{}
Let $k$ be an algebraically closed field of characteristic $p>0$, and
let $G$ be a connected, simply connected semisimple algebraic
$k$-group.  The non-0 weights of a maximal torus $T \le G$ on
$\mathfrak{g} = \Lie(G)$ form an indecomposable root system $R$ of
rank $r = \dim T$ in the character group $X = X^*(T)$.  Since $G$ is
simply connected, $X$ identifies with the full weight lattice of $R$
as in section \ref{sec:root-sys}.  We fix a choice of simple roots $S$
and positive roots $R^+$. The dominant weights are denoted $X^+$. The
group $G$ is assumed to be \emph{quasisimple}; i.e. the root
system $R$  is indecomposable.

\subsection{}
\label{sub:simple-G-mods}
For each dominant weight $\lambda \in X^+$, the space of global
sections of the corresponding line bundle on the flag variety affords
an indecomposable rational $G$-module $H^0(\lambda)$ with simple
socle.  The modules $L(\lambda) = \soc H^0(\lambda)$ comprise
all of the simple rational modules for $G$ (and are pairwise
non-isomorphic).

The character of each $H^0(\lambda)$ is the same as in characteristic 0; 
hence in particular $\dim_k H^0(\lambda)$ is given by
the Weyl degree formula, whose value at $\lambda$ we denote $d(\lambda)$ as in
\ref{sub:Weyl-degree-estimates}.

\subsection{}
\label{sub:untwisting}

Any dominant $\lambda$ may be written as a finite sum $\sum_{i \ge 0}
p^i \lambda_i$ with each $\lambda_i$ a \emph{restricted} weight.
Recall that a dominant weight $\mu$ if $\langle \mu,\alpha^\vee
\rangle < p$ for all simple roots $\alpha$.  Steinberg's tensor
product theorem says:
\begin{equation*}
  L(\lambda) \iso L(\lambda_0) \tensor L(\lambda_1)^{[1]} 
  \tensor L(\lambda_2)^{[2]} 
  \tensor \cdots
\end{equation*}
where for a $G$-module $V$, $V^{[m]}$ standards for the $m$-th
Frobenius twist of $V$.

For $d \ge 1$, let $G_d$ be the $d$-th Frobenius kernel of $G$.  Let
$V$ be a rational $G$-module and $m \ge 1$. If there is a rational $G$
module $W$ with $W^{[m]} \iso V$, we regard $W$ as the Frobenius
\emph{untwist} $W = V^{[-m]}$ of $V$. Now regard $V$ as a module for
$G_d$.  Since $G_d$ is a normal subgroup scheme, $G$ acts on
$V^{G_d}$; since $G_d$ acts trivially on this $G$-module, there is an
untwisted rational $G$-module $(V^{G_d})^{[-d]}$.  It follows that
there is an untwist $H^i(G_d,V)^{[-d]}$ for all $i \ge 0$.

Consider now two $G$-modules $V_1$ and $V_2$, and form $W = V_1
\tensor V_2^{[d]}$. The Frobenius kernel $G_d$ acts trivially on
$V_2^{[d]}$, so that 
\begin{equation}
  \label{eq:Frob-coho-untwist}
  H^i(G_d,W)^{[-d]} \iso H^i(G_d,V_1)^{[-d]}\tensor V_2
\end{equation}
as $G$-modules for every $i \ge 0$. 

\subsection{}
\label{sub:general-G-cohomology}
Let $W_p \le \widehat W_p$ be as in section \ref{sec:root-sys} (for $l
= p$), let $C = C_p \cap X^+$ denote the dominant weights in the
lowest alcove, and let $\bar C = \bar C_p \cap X^+$ ($\bar C_p$ is the
closure in $X_\R$).

\begin{prop}
  Let $\lambda \in X^+$.
  \begin{enumerate}
  \item[(a)]   If $H^i(G,L(\lambda)) \not = 0$ for some $i \ge 0$. then
  $\lambda \in W_p \da 0$.
  \item[(b)]   If $H^i(G_1,L(\lambda)) \not = 0$ for some $i \ge 0$,
    then $\lambda \in \widehat{W_p} \da 0$.
    \item[(c)] $H^i(G,H^0(\lambda)) = 0$ for all $i > 0$.
    \item[(d)] If $\lambda \in \bar C$, then $L(\lambda) =
      H^0(\lambda)$; in particular, $\dim L(\lambda) = d(\lambda)$.
  \end{enumerate}
\end{prop}

\begin{proof}
  (a) follows from the \emph{linkage principle} for $G$ \cite[Cor.
  II.6.17]{JRAG}, and (b) from the linkage principle for $G_1$
  \cite[Lemma II.9.16]{JRAG}. (c) follows from \cite[II.4.12]{JRAG}.
  (d) follows from \cite[II.6.13,II.5.10]{JRAG}.
\end{proof}

\section{The Lie algebra and the cohomology of $G_1$}

We want to describe explicitly the cohomology $H^*(G_1,k)$ in degree
$\le 2$.  For this, we need some information on the Lie algebra
$\glie$.

\subsection{}
\label{sec:lie-alg-description}

Recall that the prime $p$ is \emph{bad}[=not good] for the
indecomposable root system $R$ if one of the following holds: $p=2$
and $R$ is not of type $A_r$; $p=3$ and $R$ is of type $G_2$,$F_4$, or
$E_r$; $p=5$ and $R$ is of type $E_8$.

The prime $p$ is \emph{very good} if it is not bad, and in case
$R=A_r$, if also $p$ does not divide $r+1$.

Application of the summary in
\cite[0.13]{hum-conjugacy} yields:

\begin{leml}
  Assume that $p$ is \emph{very good.}  Then $\glie$ is a simple Lie
  algebra. The adjoint $G$-module is simple, self-dual, and isomorphic
  with $L(\tilde \alpha)$ where $\tilde \alpha$ is the dominant long
  root.
\end{leml}  
Notice that if $p>h$, then $p$ is very good.

\begin{leml}
  Assume that $p\ge h$. If $W$ is any
  $G$-module, then $\Hom_G(\glie,W^{[d]}) = 0$ for $d \ge 1$.
\end{leml}

\begin{proof}
  When $p >h$ this follows since $\glie$ is a simple $\glie$-module
  with restricted highest weight.  When $p=h$, we have $R=A_{p-1}$.
  Since $G$ is simply connected, we have $\glie = \lie{sl}_p$. Thus
  $\glie$ is an indecomposable $G$-module with unique simple quotient
  $L(\tilde \alpha)$, and the lemma follows.
\end{proof}

\subsection{}
\label{sec:null-cone}
Let $B$ be a Borel subgroup of $G$, and let $\ulie$ be the nilradical
of $\Lie(B)$.  Regarding $\ulie^*$ as a $B$-module, we get a vector
bundle on $G/B$ which we also write as $\ulie^*$. According to
\cite[3.8]{AJ1}, the formal character of the $G$-module
$H^0(G/B,\ulie^*)$ is $\chi(\tilde \alpha) = \text{ch}(\glie^*)$.  

Let $\NN \subset \glie$ be the nilpotent cone. There is by
\cite[3.9]{AJ1} an injective homomorphism of graded algebras $k[\NN]
\to H^0(G/B,S\ulie^*)$
\begin{lem}
  For simply connected, quasisimple algebraic groups $G$, $\glie^*
  \iso k[\NN]_1 \iso H^0(G/B,\ulie^*)$.
\end{lem}

\begin{proof}
  Let $I(\NN) \normal k[\glie] = S\glie^*$ be the (homogeneous)
  defining ideal of the variety $\NN$.  We need to show that $I(\NN)_1
  = 0$. If not, then $\NN \subset V \subset \glie$ for some proper
  $G$-submodule $V$.  A look at the summary in
  \cite[0.13]{hum-conjugacy} shows that, since $G$ is simply
  connected, the only $G$-submodules of $\glie$ have dimension 0 or 1.
  On the other hand, by \cite[Theorem 6.19]{hum-conjugacy}, the
  variety $\NN$ has codimension $\text{rank}(G)$ in $\glie$ and so
  clearly can't be contained in a 1 dimensional linear subspace!
\end{proof}

\begin{rems}
  \begin{enumerate}
  \item Here is a fancier result which implies the lemma if we assume
    that the prime $p$ is good for $G$. Since $G$ is simply connected
    and $p$ is good, the Springer resolution
    $$\varphi:\tilde\NN = G \times^B \ulie \to \NN$$
    given by
    $(g,X) \mapsto \Ad(g)(X)$ is a \emph{desingularization}, hence in
    particular a birational map; see \cite[Theorem 6.3 and Theorem
    6.20]{hum-conjugacy}.  Again since $G$ is simply connected and $p$
    is good, the variety $\NN$ is normal (\cite[Theorem
    4.24]{hum-conjugacy}).  Standard arguments then yield an
    isomorphism of graded algebras $k[\NN]
    \xrightarrow[\iso]{\varphi^*} \Gamma(\tilde\NN,\mathcal{O}_{\tilde
      \NN})$.  Finally, the projection $\tilde \NN \to G/B$ is an
    affine morphism, so that $\Gamma(\tilde\NN,\mathcal{O}_{\tilde
      \NN}) = H^0(G/B,S\ulie^*)$ as a graded algebra.
  \item On the other hand, if $G = PGL_r$, and $p \vert r$, one can
    find a linear form on $\glie$ that vanishes on $\NN$, hence there
    can be no isomorphism $k[\NN]_1 \to H^0(G/B,\ulie^*)$ (compare
    formal characters). So the lemma can fail when $G$ is not simply
    connected. [Note that $\varphi$ is not birational in this example.
    One can show that there is a $G_{sc}$-isomorphism $\psi:\tilde
    \NN_{sc} \to \tilde \NN$ (using some obvious notations). We get
    therefore a commuting diagram:
    \begin{equation*}
      \xymatrix{
      \tilde \NN \ar[r]^{\varphi_{sc} \circ \psi^{-1}} \ar[dr]_\varphi 
      &  \NN_{sc} \ar[d]^\gamma \\
      & \NN
      }
    \end{equation*}
    The map $\varphi_{sc} \circ \psi^{-1}$ is birational. Since
    $\gamma^*k(\NN) \subset k(\NN_{sc})$ is a proper purely
    inseparable extension, so too is $\varphi^*k(\NN) \subset
    k(\tilde \NN)$.]
  \end{enumerate}
\end{rems}

\begin{prop}
  \begin{enumerate}
  \item If $p \not = 2$ or if $R$ is not of type $C_r$, then
    $H^1(G_1,k) = 0$.
  \item Assume that $p \ge h$. Then $H^2(G_1,k)^{[-1]} \iso \glie^*$ as
    $G$-modules.
  \end{enumerate}
\end{prop}

\begin{proof}
  For (1) see \cite[Lemma II.12.1]{JRAG}. For (2), first suppose that
  $p>h$.  By \cite[3.7,3.9]{AJ1}, there is a $G$-equivariant
  isomorphism of graded rings $k[\NN]' \iso H^*(G_1,k)^{[-1]}$ where
  $k[\NN]'$ is again the graded coordinate ring of $\NN$, but with the
  linear functions on $\glie$ given degree 2. The claim now follows
  from the lemma.
  
  When $p=h$, apply \cite[Cor. 6.3]{AJ1} to see that
  $H^2(G_1,k)^{[-1]} \iso H^0(G/B,\ulie^*)$; the claim follows 
  again from the lemma in this case.
\end{proof}

\section{Low dimensional modules for $G$}

\subsection{}
\label{sub:small-generalities}

We recall first some facts about low dimensional modules established
in \cite{J3a} and \cite{serre}.
\begin{prop}
  Let $L$ be a simple non-trivial restricted $G$ module with highest
  weight $\lambda$. Suppose that $\dim L \le p$.
  \begin{enumerate}
  \item[(a)] $\lambda \in \bar C$.
  \item[(b)] $\lambda \in C$ if and only if $\dim_k L < p$.
  \item[(c)] $h \le p$. If moreover $\dim L < p$, then $h < p$.
  \item[(d)] If $R$ is not of type $A$ and $\dim L = p$, then $h<p$.
    If $p=h$ and $\dim L = p$, then $R = A_{p-1}$ and $\lambda = \varpi_i$ with
    $i \in \{1,p-1\}$.
  \end{enumerate}
\end{prop}
\begin{proof}
  (a) follows from \cite[Lemma 1.4]{J3a}, and (b) from
  \cite[1.6]{J3a}, see also \cite{serre}.  For (c), note first that
  (a) implies $\dim L = d(\lambda)$ by Proposition
  \ref{sub:general-G-cohomology}(d).  If $\lambda \in \bar C \setminus
  C$, then (a) and (b) imply that $\dim L = p$, whence $p=h$ follows from
  Weyl's degree formula. (c) now follows since $C$ is empty if $p<h$
  and $C = \{0\}$ if $p=h$.
  
  In \cite[1.6]{J3a}, Jantzen made a list of all simple restricted
  modules for $G$ with dimension $p$.  Inspecting that list
  yields (d).
\end{proof}

\subsection{Vanishing results when $\glie$ acts non-trivially}
\label{sub:G-vanish-for-G1-non-triv}
Let $L$ be a simple module for $G$.
\begin{prop}
  If $G_1$ (equivalently: $\glie$) acts non-trivially on $L$ and $\dim
  L \le p$, then $H^i(G,L) = 0$ for all $i \ge 0$.
\end{prop}

\begin{proof}
  Write the highest weight of $L$ as $\lambda = \mu_1 + p\mu_2$ with
  $\mu_1$ restricted.  Since $L^\glie = 0$, we have $\mu_1 \not = 0$.
  Since $p \ge \dim L \ge \dim L(\mu_1)$, Proposition
  \ref{sub:small-generalities} implies that $\mu_1 \in \bar C$ and
  that $h \le p$. We have in particular that $L(\mu_1) = H^0(\mu_1)$,
  hence the proposition will follow from Proposition
  \ref{sub:general-G-cohomology} if we show that $\mu_2$ is 0.
   
  If $\dim L = p$, Steinberg's tensor product theorem gives $\mu_2 =
  0$.  If $\dim L < p$ then \ref{sub:small-generalities} shows that
  $p< h$ and $\mu_1 \in C$. If $H^i(G,L) \not = 0$ for some $i$, then
  $\lambda \in W_p \da 0$ by the linkage principle, whence $\mu_1 \in
  W \da 0 + pX = \widehat{W_p}\da 0$.  Now Proposition
  \ref{sub:Weyl-degree-exceptions} applies; it shows that $\dim
  L(\mu_1) = p-1$ whence we have $\mu_2 = 0$ by another application of
  Steinberg's theorem.
\end{proof}

\subsection{Second cohomology}
\label{sub:H2-tool}

Here we prove our main tool for describing second cohomology;
first we require the following:
\begin{lem}
  Let $E_2^{p,q} \implies H^{p+q}$ be a convergent, first quadrant
  spectral sequence.
  \begin{enumerate}
  \item If $E_2^{0,1} = E_2^{1,1} = E_2^{0,2} = 0$, then $H^2 \iso E_2^{2,0}$
  \item If $E_2^{1,0} = E_2^{1,1} = E_2^{2,0} = 0$, then $H^2 \iso E_2^{0,2}$.
  \end{enumerate}
\end{lem}

\begin{proof}
  We verify (1), the argument for (2) is the same.  We must show that
  $E_\infty^{2,0} \iso E_2^{2,0}$; first note that $E_3^{2,0}$ is the
  cohomology of the sequence
  \begin{equation*}
    E_2^{0,1} \to E_2^{2,0} \to E_2^{4,-1}
  \end{equation*}
  from which we get $E_3^{2,0} \iso E_2^{2,0}$. For any first quadrant
  spectral sequence one has (by similar reasoning) that $E_a^{2,0}
  \iso E_{a+1}^{2,0}$ for $a > 2$, so we get the desired isomorphism.
\end{proof}

\begin{theorem}
  Suppose that $p \ge h$. Let $V$ be a $G$-module for which $H^i(G,V)
  = 0$ for $i=1,2$, and let $d \ge 1$.
  \begin{enumerate}
  \item $H^1(G,V^{[d]}) = 0$, and
  \item $H^2(G,V^{[d]}) \iso \Hom_G(\glie,V)$.
  \end{enumerate}
\end{theorem}

\begin{proof}
  The Frobenius kernel $G_1$ is a normal subgroup of $G$; thus there
  is a Lyndon-Hochschild-Serre spectral sequence computing
  $H^i(G,V^{[d]})$ which in view of
  \ref{sub:untwisting}~\eqref{eq:Frob-coho-untwist} has the form
  \begin{equation*}
    E_2^{s,t} = H^s(G,H^t(G_1,V^{[d]})^{[-1]}) = 
    H^s(G,H^t(G_1,k)^{[-1]} \tensor V^{[d-1]})  
  \end{equation*}
  If $t = 1$, $E_2^{s,t} = 0$ by Lemma \ref{sec:null-cone}(1).

  There is an exact sequence of the form \cite[I.4.1(4)]{JRAG}
  \begin{equation*}
    0 \to E_2^{1,0} \to H^1(G,V^{[d]}) \to E_2^{0,1}=0.
  \end{equation*}
  Thus $H^1(G,V^{[d]}) \iso E_2^{1,0} \iso H^1(G,V^{[d-1]})$.  We get
  now (1) by induction on $d$.
  
  Lemma \ref{sec:null-cone}(2) shows now that $H^2(G_1,k) \iso
  \glie^*$. Thus, the only possible non-0 $E_2$ terms of total degree 2 are
  \begin{align*}
    E_2^{0,2} =& H^0(G,\glie^* \tensor V^{[d-1]}) = \Hom_G(\glie,V^{[d-1]}) \\
    E_2^{2,0} =& H^2(G,V^{[d-1]}).
  \end{align*}
  
  For $d>1$, we apply \ref{sec:lie-alg-description} Lemma B to see
  that $E_2^{0,2} = 0$ whence $H^2(G,V^{[d]}) \iso E_2^{2,0} =
  H^2(G,V^{[d-1]})$ by part (1) of the lemma; thus (2) will follow
  provided it holds for $d=1$. In that case, we have $E_2^{2,0} = 0$
  by assumption, and the result just proved in part (1) shows that
  $E_2^{1,0}=0$.  Thus part (2) of the lemma applies; it shows that
  $H^2(G,V^{[1]}) \iso E_2^{0,2} = \Hom_G(\glie,V)$ as desired.
\end{proof}

\subsection{The second cohomology of small modules.}
\label{sub:second-vanishing}

Let $L = L(\lambda)$ be a simple $G$-module, and suppose that $\dim L
\le p$.  Proposition \ref{sub:G-vanish-for-G1-non-triv} showed that
the vanishing of cohomology for $L$ is a consequence of the linkage
principle when $\lambda \not \in pX$. However, if $\lambda \in p\Z R$,
$\lambda$ is linked to 0, so the linkage principle does not yield
vanishing. The following result shows that, despite the linkage of
$\lambda$ and $0$ in this case, the second cohomology is usually 0.
\begin{theorem}
  Let $L$ be a simple $G$-module with $\dim L \le p$. If $H^2(G,L)
  \not = 0$, then $L \iso \glie^{[d]}$ for some $d \ge 1$.
\end{theorem}

\begin{proof}
  Let $L'$ be such that $L \iso (L')^{[d]}$ for $d \ge 0$, and such that
  $\glie$ acts non-trivially on $L'$. We have by
  \ref{sub:small-generalities} that $p \ge h$.  Also, we have by
  Proposition \ref{sub:G-vanish-for-G1-non-triv} that $H^i(G,L') = 0$
  for $i \ge 1$.  If $d=0$, we are done. If $d>1$, then Theorem
  \ref{sub:H2-tool} applies, and we get that
  \begin{equation*}
    H^2(G,L) \iso \Hom_G(\glie,L').
  \end{equation*}
  
  We get by Proposition \ref{sub:small-generalities} that $p > h$
  unless $R = A_{p-1}$ and $L' = L(\varpi_i)$ with $i \in \{1,p-1\}$.
  If $p>h$, then $\glie$ is a simple $G$-module by Lemma
  \ref{sec:lie-alg-description}. So if $\Hom_G(\glie,L') \not = 0$
  then $L' \iso \glie$ whence $L \iso \glie^{[d]}$ as claimed.
  
  In the remaining case, one must just note that weight considerations
  yield $\Hom_G(\glie,L(\varpi_i)) =0$ for $i=1,p-1$, whence $H^2(G,L) =
  0$.
\end{proof}

\subsection{The second cohomology of twists of the adjoint module}
\label{sub:adjoint-second}

The first assertion of Theorem \ref{intro:lie-alg} of the introduction
follows from the following:

\begin{prop}
  Assume that $p > h$. Then $H^1(G,\glie^{[d]}) = 0$ and
  $H^2(G,\glie^{[d]}) \iso \End_G(\glie)$ has dimension 1 for $d \ge
  1$.
\end{prop}

\begin{proof}
  Since $p>h$, Lemma \ref{sec:lie-alg-description} shows that $\glie$
  is the simple module with highest weight $\tilde \alpha$. It follows
  that $\glie = H^0(\tilde \alpha)$, and thus that $H^i(G,\glie) = 0$
  for $i \ge 1$ by Proposition \ref{sub:general-G-cohomology}.  The
  proposition now follows from Theorem \ref{sub:H2-tool}.
\end{proof}

\begin{rem}
  Note that $\dim \glie > h$ (in fact, $\dim \glie = (h+1)r$ where
  $r$ is the rank of $G$).  So we get also: if $\dim \glie \le p$, then
  $\dim H^2(G,\glie^{[d]}) = 1$ for $d \ge 1$.
\end{rem}

\subsection{A second proof}
\label{sub:finite-obstruction}

Here we give a second proof of the non-vanishing of $H^2$ for twists of
the adjoint module; the result proved here verifies the remaining
assertion of Theorem \ref{intro:lie-alg} of the introduction.  We have
included the argument since it offers some ``explanation'' for the
non-vanishing.

The group $G$ arises by base change from a split reductive group
scheme $\GG$ over $\Z$.  Let $\Z_p$ be the complete ring of $p$-adic
integers, and let $\Q_p$ be its field of quotients.  For any finite
field extension $F$ of $\Q_p$, let $\oo$ denote the integers in $F$.
The residue field $\oo/\mathfrak{m}$ may be identified with the
extension $\F_q$ of $\F_p$.
  
Let $K$ denote the group of points $\GG(\oo)$ regarded as a subgroup
of $\GG(F)$.  Since $\GG$ is smooth, the reduction homomorphism $K \to
\GG(\F_q)$ is surjective (see \cite[3.4.4]{tits:reductive/local}).

For $n \ge 1$, let $K_n \subset K$ be the kernel of the map $K \to
\GG(\oo/\mathfrak{m}^n)$. Note that $K/K_1 = \GG(\F_q)$ acts by
conjugation on each quotient $K_n/K_{n+1}$.

\begin{prop}
  \label{stmt:p-adic-observation}
  \begin{enumerate}
  \item[(a)] There is for each $m \ge 1$ a canonical isomorphism
    $K_m/K_{m+1} \iso \glie_{\F_q}$ as representations for $\GG(\F_q)$,
    where $\glie_{\F_q}$ is the Lie algebra of $\GG_{\F_q}$.
  \item[(b)] If $H^2(\GG(\F_q),\glie_{\F_q}) = 0$, the exact sequence of groups
    \begin{equation*}
      1 \to K_1 \to K \to \GG(\F_q) \to 1
    \end{equation*}
    splits. 
  \item[(c)] There is a $p$-power $q_0$, depending only on the root
    system $R$ of $G$, such that $H^2(\GG(\F_q),\glie_{\F_q}) \not =
    0$ whenever $q \ge q_0$.
  \item[(d)] There is an integer $a_0 \ge 1$ such that $H^2(G,\glie^{[a]})
    \not = 0$ whenever $a \ge a_0$.
  \end{enumerate}
\end{prop}
\begin{proof}
  (a) Follows from \cite[II.\S4.3]{Demazure-Gabriel}.
  (b) Since $K_1$ is a pro-$p$ group \cite[Lemma 3.8]{Platonov}, this follows from
  \cite[Lemma 3]{serre:localclassfieldtheory}.
  
  (c) Choose a $\Q_p$ vectorspace $V$ and a non-trivial faithful
  $\Q_p$-rational representation $\GG_{\Q_p} \to \GL(V)$.  For each
  extension $F$ of $\Q_p$ with integers $\oo$, the group $K =
  \GG(\oo)$ is a subgroup of (the group of $F$-points of) $\GL(V_F)$.
  If $H^2(\GG(\F_q),\glie_{\F_q}) = 0$, the sequence in (b) is split
  and $V_F$ is a non-trivial $F[\GG(\F_q)]$-module.
  
  Since $F$ has characteristic 0, it is well known that the minimal
  dimension of a non-trivial $F[\GG(\F_q)]$ module is bounded below by
  the value $f(q)$ of a polynomial $f \in \Q[x]$, depending only on
  $G$, for which $f(q) \to \infty$ as $q \to \infty$.  We may choose
  $q_0$ such that $f(q)>\dim_{\Q_p}V$ for each $q > q_0$, and (c)
  follows at once.

  (d) now follows from (c) and \cite[Cor. 6.9]{CPSv}.
\end{proof}
 
\section{Small simple modules for $G_1$}
\label{sec:G1}

Combining results of \cite{kumar99:_froben_split} with the results
recorded in \ref{eq:small-possibilities}, we obtain some explicit
results on $G_1$ cohomology of low dimensional simple modules:
\begin{prop}
  Let $L$ be a non-trivial simple $G_1$ module with $\dim \le p$.
  Assume for some $i \ge 0$ that $H^i(G_1,L) \not = 0$. Then $\dim L
  = p-1$.  Moreover, there is a quadruple $(R,\lambda,i(0),V)$ in the
  following table for which $R$ is the root system of $G$, $\lambda$
  the high weight of $L$, $i \ge i(0)$ and $H^{i(0)}(G_1,L)^{[-1]}
  \iso V$ as $G$-modules.

  \begin{equation*}
    \begin{array}[t]{lllll}
      R & \lambda & i(0) & H^{i(0)}(G_1,L)^{[-1]}\\
      \hline \hline
      A_1 & (p-2)\varpi_1 & 1 & L(\varpi_1) \\
      A_{p-2} & \varpi_1, \varpi_{p-2} & p-2 & L(\lambda) \\
      C_{(p-1)/2} \quad \text{$p$ odd} & \varpi_1 & p-2 & L(\lambda) 
    \end{array}
  \end{equation*}
\end{prop}

\begin{proof}
  By \cite[Prop. II.3.14]{JRAG}, $L = \res^G_{G_1} L(\lambda)$ for
  some restricted dominant weight $0 \not = \lambda$. Thus
  $L(\lambda)$ is a restricted, simple $G$ module with dimension $\le
  p$.  It follows from Proposition \ref{sub:small-generalities} that
  $h \le p$, that $\lambda \in \bar C$, and that $L = H^0(\lambda)$ as
  modules for $G$.
  
  Suppose that $H^i(G_1,L) \not = 0$ for some $i$.  By the linkage
  principle for $G_1$ (Proposition \ref{sub:general-G-cohomology}(b)),
  we must have $\lambda \in \widehat{W_p} \da 0$, hence $\lambda \in
  C$.  This implies that $h < p$. Proposition \ref{sub:fund-group}
  shows that $\lambda = (p-h)\varpi_i = w_0w_i\da 0 + p\varpi_i$ for
  some $i \in J$, and Proposition \ref{prop:degree-estimates} yields
  $\dim L = p-1$ and lists the possible pairs $(R,\lambda)$.
  
  For $h < p$, Kumar, Lauritzen and Thomsen \cite[Theorem
  8]{kumar99:_froben_split} have extended a result of Andersen and
  Jantzen \cite[3.7]{AJ1}; this result implies in particular that the
  minimal degree for which $H^*(G_1,L)$ is non-0 is $\ell(w_0w_i)$,
  and that
  \begin{equation*}
    H^{\ell(w_0w_i)}(G_1,L)^{[-1]} \iso 
    H^0(\varpi_i).
  \end{equation*}
  It is straightforward to compute for each pair $(R,\lambda)$ the
  length $\ell(w_0w_i)$; one gets in this way the result.
\end{proof}

\begin{rem}
  The Theorem implies the fact (used by Jantzen in the proof of
  \cite[Lemma 1.7]{J3a}) that $H^1(G_1,L) = 0$ for all simple $G_1$
  modules $L$ with $\dim L \le p$. The argument used by Jantzen
  there relied on the calculations of $H^1$ carried out
  in \cite{J2}.
\end{rem}

\providecommand{\bysame}{\leavevmode\hbox to3em{\hrulefill}\thinspace}

\end{document}